\documentstyle{amsart}
\newtheorem{rem}{Remark}
\newtheorem{theo}{Theorem}

\newtheorem{prop}{Proposition}
\newtheorem{coro}{Corollary}
\newtheorem{exam}{Example}

\begin{document}
\numberwithin{equation}{section}
\title[Generalization of second  trace form of  central simple algebras]
{Generalization of the second  trace form of  central simple
algebras in characteristic two.}
\author{Ana-Cecilia de la Maza}
\address{Instituto de Matem\'atica y F\'{\i}sica,
 Universidad de Talca, Casilla 747, Talca, Chile\\
 Fax 56-71-200392}

\thanks{This work  was supported by Fondecyt grants N$^{o} 1010205$, and by
Programa formas extremas y representaci\'on de formas
cuadr\'aticas, Universidad de Talca.}
\begin{abstract}Let $F$ be a field with characteristic two.
We generalize the second trace form for central simple algebras
with odd degree over $F$. We determine the second trace form and
the Arf invariant and Clifford invariant for tensor products of
central simple algebras.
\end{abstract}


\maketitle
\section{Introduction}
Let $A$ be a central simple algebra over a  field  $F$. For each
$a\in A$, let $Prd_{A,a}(x)=
x^n-t_1(a)x^{n-1}+t_2(a)x^{n-2}+\cdots+(-1)^nt_n(a)$ be the
reduced polynomial of $a$ (so $t_1(a)$ is the reduced trace
$Trd_{A}(a)$ and $t_n(a)$ is the reduced norm $Nrd_A(a)$ of $a$.)
Put $q_1(x)=t_1(x^2)$ and $q_2(x)=t_2(x)$. When the characteristic
of $F$ is not equal to $2$, the trace form $(A,q_1)$ and the
second trace form $(A,q_2)$ are nonsingular quadratic forms (see
\cite{l}, \cite{lm} and \cite{t}). If the characteristic is $2$,
however, then the trace form $(A,q_1)$ has rank zero (and is
therefore singular). In this situation, the second trace form
$(A,q_2)$ is nonsingular if the degree of $A$ is even (see
\cite{bf}) but it is necessarily singular if the degree is odd.

In this article, we extend the definition of the second trace form
to the case that the degree of $A$ is odd and study the behavior
with respect to tensor products. Our definition is similar to the
way in which Revoy defined the trace form $T_{E/F}$ for a field
extension $E$ of $F$ in \cite{r}.

One reason for wanting to have the notion of a second trace form
for central simple algebras of odd degree is the following. When
we have two central simple algebras, $A $ with even degree and $B$
with odd degree, the tensor product $A\otimes B$ has even degree.
Hence, we have two nonsingular second trace forms $(A,q_2)$ and
$(A\otimes B,q_2)$. Now, if we extend the definition of the second
trace form to a nonsingular quadratic form over the central simple
algebras with odd degree, then is possible to obtain $(A\otimes B,
q_2)$ through $(A,q_2)$ and $T_{B/F}$, where $T_{B/F}$ is the
second trace form for odd component $B$. Decompositions of this
type naturally appear when decomposing central simple algebras in
terms of a matrix algebra and a division algebra, and in the
primary decomposition of division algebras. We show (cf.
Theorem~2) that our definition of the second trace form for
central simple algebras of odd degree is compatible with such
tensor product decompositions.

Let us now describe the contents of the paper in more detail. In
Section 2, we define the second trace form $T_{A/F}$ for central
simple algebras $A$ (with any parity) over a field $F$ of
characteristic two. We prove that it is a nonsingular quadratic
form over $F$. In Section 3, we calculate this form for a crossed
product. In particular, we prove that---given a Galois field
extension $E/F$ of odd degree---the second trace form for the
crossed product algebra $(E,G,\Phi)$ is Witt equivalent to the
second trace form $T_{E/F}$ of Revoy.  We conclude that section
with a relatively simple  criterion that uses the trace form in
order to  recognize in many cases that a given field extension is
not Galois. Note that in general this is not easy for fields of
characteristic two. In section 4, we study the second trace form
for tensor products of central simple algebras over $F$ (for any
parity). As an application, we determine the Arf invariant and the
Clifford invariant for tensor products in Section 5 (see Refs.
\cite{a}, \cite{b} and \cite{sah} for properties of these
invariants).

\section{Second trace form}

In this section we will define the second trace form for central
simple algebra over a field of characteristic two. We will prove
that this is a nonsingular quadratic form. First, we give some
properties and notations.

Given a central simple algebra $A$ over $F$,  the degree of $A$,
deg$A$ for short, is the integer $n$ such that $\dim_F A=n^2$
\cite[p. 236]{p}. A splitting field for $A$ is a field $E$
containing $F$  such that $A\otimes_F E$ is isomorphic to the
matrix algebra $M_n(E)$,  and a splitting representation is a
$F$-algebra isomorphism $\phi:A \rightarrow M_n(E)$. For each
$a\in A$  the reduced polynomial $Prd_{A,a}(x)$ is defined as
\begin{equation}\label{pol}
Prd_{A,a}(x):=det(xI_n-\phi(a))=
x^n-t_1(a)x^{n-1}+t_2(a)x^{n-2}+\cdots+(-1)^nt_n(a).
\end{equation}
This reduced polynomial $Prd_{A,a}(x)$ has coefficients in $F$
(i.e. it lies in $F[x]$) and it is independent of $E$ and $\phi$
(see \cite[p.295]{p}).

For each central simple algebra $A$ over $F$, $(A,t_2)$ is a quadratic space.
If in particular the characteristic of $F$ is  two, then it is easy to
prove that for each $a\in A$
\begin{equation}\label{t2}
t_2(a)= \sum_{1\leq i<j\leq n} \big(\phi(a)_{ii}\
\phi(a)_{jj}+\phi(a)_{ij}\ \phi(a)_{ji}\big)
\end{equation}
(where $\phi$ is a splitting representation) and, furthermore,
that the associated bilinear form $b_{t_2}$ satisfies
\begin{equation}\label{bt2}
b_{t_2}(x,y)=t_1(xy) +t_1(x)t_1(y) \quad (\mbox{for each } x,\
y\in A).
\end{equation}

Since the quadratic space $(E,t_2)$ is necessarily singular when
$\deg A$ is odd, we will define a reduced second trace form $T_{A/F}$
in the spirit of Revoy's definition in \cite{r} for quadratic forms
over extension fields (see section 3). To this end we first note that,
by Eq. \eqref{bt2} and the linearity of $t_1$, the spaces $F$ and
$A_0:= Ker \ t_1$ are mutually orthogonal. The reduced trace form
$T_{A/F}$ is now defined as
\begin{equation}\label{deft2}
T_{A/F}=
\begin{cases}
(A,t_2) & \text{   if  }\ \deg\ A \ \text{ is  even},
\\
(A_0,t_2) &\text{  if  }\ \deg\ A \  \text{  is odd}.
\end{cases}
\end{equation}

\begin{rem} Let $A$ and $B$ be  isomorphic central simple algebras over $F$.
Then, by Eqs. \eqref{t2} and \eqref{bt2}, it is clear that
$T_{A/F}$ and $T_{B/F} $ are isometric.
\end{rem}

Let us denote the quadratic form $ax^2+xy+by^2$ as $[a,b]$ and the
hyperbolic plane $[0,0]$ as $\mathbb{H}$.

\begin{prop}\label{mn}
For a matrix algebra $A=M_n(F)$ with $n>1$, the reduced trace form
$T_{A/F}$ is Witt equivalent to
\begin{equation}\label{matr}
\begin{cases}
\mathbb{H} & \text{ if } \quad n\equiv 0,\ 1,\ 2,\ 7\mod 8, \\
\left[ 1,1\right] & \text{ if } \quad n \equiv 3,\ 4,\ 5,\ 6\mod
8.\end{cases}\end{equation}
\end{prop}

\noindent Proof. We write   the canonical base for $M_n(F)$ as
$\{E_{ij}\}$ and put $e_i:=E_{ii}$. For each triple of natural
numbers $0<k_1\leq n$ and $0<k_2,k_3\leq \frac{n+3}{4}$ we define:
\begin{eqnarray}
V_{k_1}&:=& span \cup_{1\leq i< j\leq k_1 }\{ E_{ij},\ E_{ji}\},
\nonumber \\
V'_{k_2}&:= & span \cup_{0\leq i\leq k_2}\{
e_{4i+1}+e_{4i+2},\ e_{4i+2}+e_{4i+3}\}, \nonumber \\
W_{k_3}&:= & span \cup_{1\leq i\leq k_3} \left\{
\sum_{t=1}^{4i}e_t\ ,\
e_{4i}+e_{4i+1}+e_{4i+2}+e_{4i+3}\right\}.\nonumber
\end{eqnarray}
It is clear that each pair of these subspaces has trivial
intersection. It moreover follows from Eq. \eqref{bt2} that the
above bases for the subspaces $V_{k_1},\ V'_{k_2}$ and $W_{k_3}$
are in fact symplectic. Hence, using Eqs. \eqref{t2} and
\eqref{bt2} we obtain that
$(V_{k_1},t_2)=\frac{k_1(k_1-1)}{2}{\mathbb{H}}, \
(V'_{k_2},t_2)=(k_2+1)[1,1]$, and $(W_{k_3},t_2)=k_3{\mathbb{H}}.$
Furthermore by Eq. \eqref{bt2} the subspaces in question are
mutually orthogonal. Defining $S_{k_1, k_2,k_3}$ as  the direct
sum of those spaces, i.e. $S_{k_1,k_2,k_3}=V_{k_1}\oplus
V'_{k_2}\oplus W_{k_3}$, we have that
\begin{equation}\label{vk}
 (S_{k_1,k_2,k_3},t_2)=
(k_2+1)[1,1]+\Big(\frac{k_1(k_1-1)}{2}+k_3\Big)\,\mathbb{H}.
\end{equation}
Defining $W_0=\{\theta\} $, where $\theta$ is  the null matrix, we
extend  the definition of $S_{k_1,k_2,k_3}$ to the case $k_3=0$.
With this notation, we have that
\[A=\begin{cases}
S_{n,k-1,k-1}\oplus^\perp\Big\langle \sum_{i=1}^{4k} e_i, e_{4k}
\Big\rangle & \text{ if }  n=4k, \\
V_2\oplus^\perp \Big\langle \sum_{i=1}^{2}e_i,\ e_2 \Big\rangle &\text{ if } n=2,\\
S_{n,k-1,k-1}\oplus^\perp \Big \langle \sum_{i=1}^{4k}
e_i,e_{4k}+e_{4k+1}\Big \rangle \perp \Big\langle
\sum_{i=1}^{4k+2}e_i,e_{4k+2} \Big\rangle & \text{ if }\begin{cases} n=4k+2,\\
n\neq 2,\end{cases}
  \end{cases}
\]
and
\[A_0=\begin{cases}
S_{n,k-1,k-1}\oplus^\perp\Big\langle  \sum_{i=1}^{4k}e_i,\ e_{4k}+
e_{4k+1}\Big\rangle & \text{ if }  n=4k+1, \\
S_{n,k,k}& \text{ if } n=4k+3,
\end{cases}\]
where $\Big\langle x,y\Big\rangle$ denotes the space generated by
$x$ and $y$. Using Eq. \eqref{vk}, the fact that  $[0,0]=[1,0]$
(see \cite[p. 150]{sah}), and replacing $n=2m+1$ with $m=2k$ or
$m=2k+1$, we obtain that $T_{A/F}$ can be written as
\begin{equation*}
\begin{cases}
\left[\frac{n}{4}\right]\left[ 1,1\right] \perp \big(
2m^2-\frac{m}{2}\big)\mathbb{H}&\text{if } \quad n=2m, \\
&\\
\left[\frac{n+1}{4}\right]\left[ 1,1\right] \perp  \big(2m^2+
\left[\frac{3m}{2}\right]\big){\mathbb{H}}    & \text{if  }\quad
n=2m+1,
\end{cases}
\end{equation*}
where  $[\frac{x}{y}]$ denotes the integer part of $\frac{x}{y}$.
Hence, using that $[1,1 ]+[1,1]=2\mathbb{H}$, the proof is
complete.  \hfill $\Box$

\begin{prop}
 For each central simple algebra $A$ over $F$ with $A\neq F$, the second
 trace form $T_{A/F}$ is a nonsingular quadratic form over $F$.
\end{prop}
\noindent Proof. Let $A$ be a central simple algebra over $F$ with
$A\neq F$. Let $E$ be a splitting field of $A$ with
$\phi:A\rightarrow M_n(E)$ the splitting representation. By
\cite[p. 238]{p},we extend $\phi$ to $E$-algebra isomorphism
$\phi:A\otimes_F E\rightarrow M_n(E)$ given by $\phi(a\otimes
e)=\phi(a)e$. Clearly $\phi(A_0\otimes_F E)= M_n(E)_0$, whence the
quadratic form $(A_0\otimes_F E,t_2)=(M_n(E)_0,t_2)$ is
nonsingular for $\deg A$ odd (by Proposition \ref{mn}). As a
consequence, $(A_0,t_2)$ is nonsingular. Similarly, if $\deg\ A $
is even then $(A\otimes_F E,t_2)$ is nonsingular (again by
Proposition \ref{mn}). Hence $(A,t_2)$ is nonsingular. \hfill
$\Box$

\begin{rem} Observe that we give here an alternative proof for the
even case, already established by Berhuy and Frings in \cite[p.
4,5]{bf}.
\end{rem}

\section{Crossed product algebra}

In this section we will compute the second trace form for a
crossed product.

Given a field extension $E/F$, we denote by $T_{E/F}$ the second
trace form  due to Revoy \cite{r}, that is
\begin{equation}\label{defrev}
T_{E/F}=
\begin{cases}(E,T_2) & \text{   if  }\ \ [E:F] \ \text{
is  even},\\
(Ker \ T_1,T_2) &\text{  if  }\ \ [E:F] \  \text{  is odd},
\end{cases}\end{equation}
where for $a \in E$, $T_1(a)$ and  $T_2(a)$ denote the
coefficients of $x^{[E:F]-1}$  and  $x^{[E:F]-2}$, respectively,
in the characteristic polynomial of $a$ (so $T_1(a)$ is the trace
$tr_{E/F}(a)$ of $a$). Note that there are  two different
definitions for the second trace forms
 in the literature. The second trace form,  due
to Berg\'e and Martinet \cite{bm}, increases the dimension of the
space by $1$ using the \'etale $F$-algebra. In \cite{m} we proved
that the two definitions are Witt equivalent.

\begin{prop} Let $E/F$ be a Galois field extension, with
$Gal(E/F)=G$. Let $A=(E,G,\Phi)=\sum_{\sigma\in G} u_\sigma E$ be
the crossed product, where $\Phi$ is normalized and for each
$\sigma, \tau\in G$ and $c\in E$
\begin{equation}\label{form}
u_\sigma^{-1}cu_\sigma=\sigma(c) \mbox{ and }
\Phi(\sigma,\tau)=u_{\sigma\tau}^{-1}u_\sigma u_\tau\in
E.\end{equation} Then
\begin{itemize}
\item[i)] $(E,t_2)=(E,T_2)$, where $ T_2$ as above and
$t_2$ as in \eqref{pol};
\item[ii)] $t_1(\sum_{\sigma\in G} u_\sigma
c_\sigma)=tr_{E/F}(c_{id})$, where $c_{id}$ is the coefficient
that correspond to $id$ (the identity);
\item[iii)] if $\sigma\neq id$, then for each  $c\in E$, $u_\sigma c\in
A_0$;
\item[iv)] if $\sigma \tau\neq id$, then for each $c,\ d\in E$,
$b_{t_2}(u_\sigma c,\ u_\tau d)=0$;
\item[v)] the subspaces $E$ and  $ \langle u_\sigma c\mid  c\in E,
 \ \sigma\neq id\rangle$ are mutually orthogonal;
\item[vi)] if $\rho^2\neq id$, then for each $c\in E$,
$t_2(u_\rho c)=0$.
\end{itemize}

\end{prop}

\indent Proof. For i) see \cite[p. 297]{p}. There is  a good
splitting representation $\phi$ of $A=(E,G,\Phi)$ in \cite[p.
298]{p}, given by
\begin{equation}\label{rep}
\phi\Big(\sum_{\rho\in G}u_\rho c_\rho\Big)=[d_{\sigma \tau}]
\quad\mbox{ where } d_{\sigma \tau}=\Phi(\sigma
\tau^{-1},\tau)\,c_{\sigma \tau^{-1}}^\tau,
\end{equation}
where  for $\tau\in G$ and $c\in E$, the notation $c^\tau$
corresponds to  $\tau(c)$. Using  that $t_1$ is the trace of the
matrix $[d_{\sigma \tau}]$, we have
\begin{equation}\label{t11}
t_1\Big(\sum_{\rho\in G}u_\rho c_\rho\Big)= \sum_{\sigma= \tau}
d_{\sigma \tau}=\sum_{\sigma\in G}c_{id}^\sigma=tr_{E/F}( c_{id}).
\end{equation}
Hence ii) and iii) are true. Furthermore, using Eq. \eqref{bt2}
together with Eq. \eqref{form} we obtain iv) (and v) as a
particular case). Now, by Eqs. \eqref{rep} and \eqref{t2}, we have
\begin{equation}\label{t22}
t_2\Big(\sum_{\rho\in G}u_\rho c_\rho\Big)= \sum_{\sigma \neq
\tau}( c_{id}^\sigma c_{id}^\tau+ u_\sigma ^{-1}
u_{\sigma\tau^{-1}}u_\tau c_{\sigma\tau^{-1}}^\tau  u_\tau^{-1}
u_{\tau\sigma^{-1}}u_\sigma c_{\tau\sigma^{-1}}^\sigma).
\end{equation}

Using the fact that if $\rho^2\neq id$ and $\rho=\sigma \tau^{-1}$
for some $\sigma $ and $\tau$ in $G$, then $\rho\neq \tau\sigma
^{-1}$, we obtain in Eq. \eqref{t22} that for $c\in E$,
$t_2(u_\rho c)=0$. \hfill $\Box$

The following theorem characterizes the second trace form
$T_{A/F}$ for a crossed product $(E,Gal(E/F),\Phi)$ in terms of
the second trace form $T_{E/F}$ of Revoy.

\begin{theo} Let $E/F$ be a Galois extension with $G=Gal(E/F)$.
Let $A$ be the crossed product $(E,G,\Phi)=\sum_{\sigma\in G}
u_\sigma E$.  Then $T_{A/F}$ is Witt equivalent to
\[\begin{cases}
T_{E/F}\perp (B,t_2) & \text{ if\  } [E:F] \text{ is even}, \\
T_{E/F} & \text{ if\  } [E:F] \text{ is odd},
\end{cases}\]
where $B$ is the subspace
\[ B:=\Big\langle u_\rho e \mid e\in E,
\ \rho \in G,\ \rho^2=id, \rho\neq id\Big\rangle.\]
\end{theo}

\noindent Proof.  We can suppose that  $\Phi$ is normalized
\cite[p. 252 ]{p}. Let us write $G$ as
\[G=\{\rho_1=id,\rho_2,\cdots,\rho_s,
\sigma_1,\sigma_2,\cdots
\sigma_t,\sigma_1^{-1},\sigma_2^{-1},\cdots \sigma_t^{-1} \},\]
where $s+2t=[E:F]$ and for $i\leq j\leq s,\ \rho_j^2=id$, and for
$1\leq i\leq t,\ \sigma_i^2\neq id$. Let $B':=\langle u_{\sigma_i}
e \mid e\in E, \ 1\leq i\leq t \rangle$. It is clear that the sum
$E+B'$ is direct (when $[E:F]$ is even, $E+B'+B$ is direct).
Furthermore, by Proposition 3, $E,\ B'$ and $ B$ are mutually
orthogonal, and $B'$ is a totally isotropic subspace (with
dimension $t [E:F]$). Hence, by combining the
fact that $\dim_F B=(s-1)[E:F]$ and that $[E:F]=s+2t$, we see that:\\
- if $[E:F]$ is even, then $T_{A/F} =T_{E/F}\perp (B,t_2) \perp
t[E:F]\,\mathbb{H}$ (see \cite[p. 17]{b}),\\
- if $[E:F]$ is odd, then $s=1$, whence $T_{A/F}=T_{E/F}\perp
t[E:F]\, \mathbb{H}$. \hfill $\Box$ \vspace{0.1in}

\begin{coro} Let $E/F$ be a Galois extensions with degree $n=2m+1$.
Let $A$ be the crossed product $(E,Gal(E/F),\Phi)$. Then
\[T_{A/F}=\begin{cases}
m\,\mathbb{H} & \text{ if }\quad n\equiv 1,\ 7\mod 8, \\
\left[1,1\right] +(m-1)\,{\mathbb{H}}  & \text{ if }\quad n\equiv
3,\ 5 \mod 8.
\end{cases}\]
\end{coro}

\noindent Proof. The algebra $B=End_F(E)\cong M_n(F)$ is the
crossed product $(E,G,I)$ \cite[p. 252]{p}. Hence, by Theorem 1,
$T_{A/F}=T_{E/F}=T_{B/F}$. Using Eq. \eqref{matr}, we obtain the
result. \hfill $\Box$
\vspace{0.1in}

 In general it is not easy to
decide whether a given field extension is Galois. The following
by-product of our work on the trace form permits an elegant
partial solution to this problem.

\begin{coro}
Let $E/F$  be a field  extension with odd degree $n$.  If
$T_{E/F}$ is not  Witt equivalent to $\mathbb{H}$ when $n\equiv
1,7 \mod 8$, or not equivalent  to $[1,1]$ when $n\equiv 3,5
\mod 8$, then the extension $E/F$ is not Galois.
\end{coro}
\noindent Proof. Immediate, by Corollary 1. \hfill $\Box$

\begin{exam} Let $F={\mathbb{F}}_2(a)$ and $E=F(b)$, where
$a^2+a+1=0$ and $b^3+b+a=0$. It is easy to prove that
$T_{E/F}=[1,a]$. Note that $[1,1]\neq[1,a]$, because $1\in \wp
(F)=\{x^2+x| x\in F\}$ and $a\notin  \wp (F)$. Hence, $E/F$ is not
Galois by corollary 2.
\end{exam}

\section{Second trace form of a tensor product.}

Let $C.S(F)$ be the set of central simple algebras over $F$. For
each $A\in C.S(F)$ and any $a\in A$, we remember that the reduced
trace of $a$ is  $Trd_A(a):=t_1(a)$.

\begin{prop} Let $A$ and $B$ be  central simple algebras over $F$.
Then
\begin{itemize}
\item[i)] $1_A\otimes_F B_0$ and $ A_0\otimes_F 1_B$ are mutually
orthogonal (as quadratic subspaces of $(A\otimes_FB,t_2)$),
\item[ii)] for each $a, \ a' \in A$ and $b,\ b'\in B$ we have that
\begin{eqnarray}
Trd_{A\otimes_F B}(a\otimes b)&=& Trd_{A}(a)Trd_{B}(b), \label{trab}\\
t_2(a\otimes b)&=&(Trd_{A}(a))^2t_2(b)+(Trd_{B}(b))^2t_2(a),\label{t2ab}\\
b_{t_2}(a\otimes b,a'\otimes b')&=&
Trd_{A}(aa')b_{t_2}(b,b')+Trd_{B}(b)Trd_{B}(b')b_{t_2}(a,a')\nonumber \\
&=&Trd_{B}(bb')b_{t_2}(a,a')+Trd_{A}(a)Trd_{A}(a')b_{t_2}(b,b'),
\label{btab}
\end{eqnarray}
\item[iii)] if $\{a,\ a'\}\cap A_0\neq \phi$ and $\{b,\
b'\}\cap B_0\neq \phi$, then
\begin{equation}\label{bt}
 b_{t_2}(a\otimes b,a'\otimes b')=
b_{t_2}(a,a')b_{t_2}(b,b').
\end{equation}
\end{itemize}
\end{prop}
\noindent Proof. Let $A,\ B\in C.S(F)$, and let $E$ and $L$ be
splitting fields for $A$ and $B$, respectively, with splitting
representations $\phi:A\rightarrow M_m(E)$ and
$\varphi:B\rightarrow M_n(L)$. Clearly $EL$ is a splitting field
for $A\otimes_FB$ ($EL$  is the minimal field that contains $E$
and $L$). A corresponding splitting representation is given by
$\Phi:A\otimes_F B\rightarrow M_{mn}(EL)$ of the form
$\Phi(a\otimes b)_{ij}=\phi(a)_{kl}\,\varphi(b)_{st}$, where
$i=(k-1)n+s$ and $j=(l-1)n+t$, with $1\leq k,l\leq m+1$ and $0\leq
t\leq n-1$. Hence, by using the definition of the reduced trace
and Eqs. \eqref{t2}, \eqref{bt2}, we obtain the equations that
appear in ii). Furthermore, by Eq. \eqref{btab} and the fact that
$F$ and $Ker \ t_1$ are mutually orthogonal, we obtain i). In
order to obtain iii), we use that
$Trd_{A}(a)Trd_{A}(a')=0=Trd_{B}(b)Trd_{B}(b')$. \hfill $\Box$

The following theorem provides the second trace form of a tensor
product in terms of the second trace forms of its constituents.

\begin{theo} Let $A_1,\ A_2\in C.S(F)$, with $\deg\ A_i=n_i$ for $i=1,2$.
$T_{A_1\otimes A_2/F}$ can be represented as
\[\begin{cases}
T_{A_1/F}+ T_{A_2/F}+\frac{(n_1^2-1)(n_2^2-1)}{2}\mathbb{H}&
\text{ if \ \ \,} n_1, \ n_2
\text{ are  odd}, \\
\left[1,1\right]+\left(\frac{n_1^2 n_2^2}{2}-1\right)\mathbb{H} &
\text{ if \ \ \,} n_1,\ n_2\equiv 2\mod  4,\\
\frac{n_1^2n_2^2 }{2}\mathbb{H}& \text{ if }\quad
 n_i\equiv 0\mod 4 \text{ and } n_j \text{ is even }, \\
 T_{A_i/F}+\frac{n_i^2(n_j^2-1)}{2}\mathbb{H}& \text{ if }
  \begin{cases}
    n_i\equiv 0\mod  4 \text{ and } n_{j}  \text{ is odd}, \text{ or}\\
    n_i\equiv 2\mod  4 \text{ and } n_{j}\equiv 1\mod  4,
\end{cases} \\
\left[1,1\right]+T_{A_i/F}+\left(\frac{n_i^2(n_j^2-1)}{2}-1\right)
\mathbb{H} & \text{ if \ \ \,} n_i\equiv 2\mod  4 \text{ and }
n_{j}\equiv 3\mod  4,  \\
\end{cases}\]
where $\{i,j\}=\{1,2\}$.
\end{theo}

\indent Proof. For $j=1,2$, let
\[
\begin{array}{l}
\{e^{(j)}_i,\ f^{(j)}_i\}_{i\in I^{(j)}}, \text{ be a symplectic
basis of } (A_j)_0 \ \text{ if }
 n_j \text{ is odd}, \text{ and}\\
\{1_{A_j},\ f^{(j)}\}\cup \{e^{(j)}_i,\ f^{(j)}_i \}_{i\in I^{(j)}
} \text{ be a symplectic basis of }A_j \ \text{ if } n_j \text{ is
even}.\end{array}
\] Put
\begin{equation}\label{w}
W=\bigoplus^\perp_{i\in I^{(1)}, \ j\in I^{(2)} }
\left(\Big\langle e^{(1)}_i\otimes e^{(2)}_j, f^{(1)}_i\otimes f^{(2)}_j\Big\rangle\perp \Big\langle
e^{(1)}_i\otimes f^{(2)}_j,f^{(1)}_i\otimes e^{(2)}_j\Big\rangle\right).\end{equation}

In view of  \eqref{t2ab}, $(W,t_2)$ is hyperbolic.

Using Eqs. \eqref{trab}, \eqref{btab} and \eqref{bt} we obtain
decompositions of $(A_1\oplus A_2)_0$ and $(A_1\oplus A_2)$,
respectively in the following cases:
\begin{description}
\item[- $A_1$ and $A_2$ have odd degree]

\[(A_1\otimes A_2)_0=
\bigoplus^\perp_{i\in I^{(1)}}\Big\langle e^{(1)}_i\otimes
1_{A_2}, f^{(1)}_i \otimes 1_{A_2} \Big\rangle\perp \bigoplus^\perp_{j\in
I^{(2)}}\Big\langle 1_{A_1}\otimes e^{(2)}_j, 1_{A_1}\otimes
f^{(2)}_j\Big\rangle\perp W.\]

\item[- $A_1$ and $A_2$ have even degree]
\[A_1\otimes A_2=
\Big\langle 1_{A_1}\otimes f^{(2)}, f^{(1)}\otimes
1_{A_2}\Big\rangle\perp W\perp V' ,\]
where $\Big\langle 1_{A_1}\otimes 1_{A_2},\ 1_{A_1}\otimes e^{(2)}_j,\ 1_{A_1}\otimes f^{(2)}_j,
e^{(1)}_i\otimes 1_{A_2},\ f^{(1)}_i\otimes 1_{A_2}\Big\rangle_{\{i\in I^{(1)},\ j\in I^{(2)}\}}$
is a totally isotropic subspace of $V'$ with dimension $\dim_F
A_1+\dim_F A_2-3=\frac{1}{2} \dim_F V'$. Hence $(V', t_2)=(\dim_F
A_1+\dim_F A_2-3)\mathbb{H}$.

\item[- $A_1$ and $A_2$ have different parity] Suppose that the degree of $A_1$ is even. Then $A_1\otimes A_2$ can be written as
\[ \Big\langle 1_{A_1}\otimes 1_{A_2},f^{(1)}\otimes 1_{A_2}\Big\rangle\
\perp \bigoplus^\perp_{\{i\in I^{(1)}\}}\Big\langle e^{(1)}_i\otimes 1_{A_2}, f^{(1)}_i
\otimes 1_{A_2} \Big\rangle\perp W\perp
V'', \]
where
$\Big\langle 1_{A_1} \otimes e^{(2)}_j,1_{A_1}\otimes f^{(2)}_j\Big\rangle_{j\in J^{(2)}}$
is a totally isotropic subspace of $V''$ with dimension $(\dim_F
A_2-1)=\frac{1}{2} \dim_F V''$. Hence, $(V'', t_2)=(\dim_F
A_2-1)\mathbb{H}$.
\end{description}

Finally,  we obtain the claim using  Eq. \eqref{t2ab} and the fact
that for each $A\in C.S(F)$, $t_2(1_A)$ is given by
\[ t_2(1_A)=
  \begin{cases}
    0 & \text{if } \deg A \equiv 0, \ 1\mod  4, \\
    1 & \text{if } \deg A \equiv 2,\ 3\mod  4.
  \end{cases}\]
\hfill $\Box$

\begin{coro} Let $A$ be a central simple algebra with even
degree. If the second trace form $T_{A/F}$ is not Witt equivalence
to $[1,1]$ or to $\mathbb{H}$, then $A=M_k(D)$ with $k$ odd and
$D$ a division algebra.
\end{coro}
\noindent {\it Proof.} Since for each  $A\in C.S(F)$,
$A\simeq M_k(D)\simeq M_k(F)\otimes D$ for some division
algebra $D$, by Theorem 2 we obtain that $T_{A/F}$ is Witt
equivalent to $T_{D/F},\ [1,1]+T_{D/F},\ \mathbb{H}$ or $[1,1]$,
depending on $k$ and $\deg D$. Now, if $k$ is even then $T_{A/F}$
is Witt equivalent to $\mathbb{H}$ or to $[1,1]$. \hfill $\Box$

\begin{coro} Let $A$ be a central simple algebra over $F$.
\[T_{A\oplus A/F}\text{ is Witt equivalent to } \begin{cases}
\left[1,1\right] &\text{ if } \deg  A\equiv 2 \mod 4,\\
{\mathbb H} &\text{ otherwise}.
\end{cases}\]
\end{coro}

\section{Invariants}
In this section we will determine the Arf invariant and the
Clifford invariant for the tensor product of algebras as an
application of  Theorem 2.

For $a\in F^*$ and $b\in F$ the quaternion algebra $
(a,b]\in C.S(F)$ is defined as the algebra
with basis $\{1,\ e,\ f,\ ef\}$ satisfying $e^2=a$, $f^2+f=b$ and
$ef+fe=e$. Now, given a non degenerated quadratic form $(V,q)$
over $F$, we can rewrite $q$ as $q= \langle a_1\rangle [1,b_1]
\perp \cdots \perp \langle a_n\rangle [1,b_n]$, with $a_i\in F^*$
and $b_i\in F$.  The Arf invariant is given by
$Arf(q):=b_1+b_2+\cdots+b_n \mod \ \wp(F)$, where
$\wp(F):=\{x^2+x | x\in F\}$, considered as additive subgroup of $F$.
The Clifford invariant $C(V,q)$ is the class of the tensor products of
quaternion algebras $(a_1,b_1]\otimes \cdots \otimes(a_n,b_n]$ in
the Brauer group $Br(F)$. Note that if $a\neq 0$ then
\[C([a,b])=C(\langle a\rangle [1,ab])=(a,ab]=((a,b))_F,\]
where $((a,b))_F$ \cite[p.25]{kmrt} is the algebra generated by
$r$ and $s$ satisfying
\[r^2=a,\ s^2= b,\ rs+sr=1 .\]

We need the following result of Berhuy and Frings \cite{bf}.

\begin{theo}[Berhuy-Frings \cite{bf}]
Let $F$ be a field of characteristic two, $n\geq 2$ an even
integer and $A\in C.S(F)$ an algebra of degree $n$ over $F$. Then
we have $Arf(T_{A/F})=[\frac{n}{4}]$ and $C(A,
T_{A/F})=[A]^{\frac{n}{2}}$, where $[\frac{n}{4}]$ denotes the
integer part of $\frac{n}{4}$ and $[A]$ denotes the class of $A$
in the Brauer group $Br(F)$.
\end{theo}

The following theorem gives the Arf invariant and the Clifford
invariant of the second trace form of a tensor product in terms of
the corresponding invariants of its constituents.

\begin{theo}

 Let $A_1,\ A_2\in C.S(F)$, with $\deg\ A_i=n_i$.  Then
\[Arf(T_{A_1\otimes A_2/F}) =
\begin{cases}
 Arf(T_{A_1/F}) + Arf(T_{A_2/F})& \text{\ if \ } \ n_1 n_2 \text{\  is odd}, \\
 \left[\frac{n_1n_2}{4}\right]& \text{otherwise},
\end{cases}\]
and
\[C(T_{A_1\otimes A_2/F})=
\begin{cases}
C(T_{A_1/F})\cdot C(T_{A_2/F})&\text{ if \ \ \,} n_1, \ n_2 \text{ are  odd}, \\
((1,1))_F & \text{ if \ \ \,} n_1,\ n_2 \equiv 2\mod 4,\\
1& \text{ if }\quad n_i\equiv 0\mod  4\text{ and } n_j \text{ is even}, \\
\left[A_i\right]^{\frac{n_i}{2}}& \text{ if }
  \begin{cases}
    n_i\equiv 0\mod  4 \text{ and } n_{j}  \text{ is odd}, \text{ or} \\
    n_i\equiv 2\mod  4 \text{ and } n_{j}\equiv 1\mod  4,
  \end{cases} \\
((1,1))_F\cdot \left[A_i\right]^{\frac{n_i}{2}}& \text{ if \ \ \,}
n_i\equiv 2\mod  4 \text{ and }
n_{j}\equiv 3\mod  4, \\
\end{cases}\] where $\{i,j\}=\{1,2\}$.
\end{theo}

\indent Proof. When $A_1$ or $A_2$ have  even degree, we have that
$A_1\otimes A_2$ also has even degree. Hence, by Theorem 3,
$Arf(T_{A_1\otimes A_2/F})=[\frac{n_1n_2}{2}]$. But, if $A_1$ and
$A_2$ have odd degree, then by Theorem 2, $T_{A_1\otimes
A_2/F}=T_{A_1/F}+T_{A_2/F}$. Hence, in this case
$Arf(T_{A_1\otimes A_2/F}) = Arf(T_{A_1/F}) + Arf(T_{A_2/F})$. In
order to obtain the Clifford invariants, we need Theorem 2, the
fact that $C({\mathbb{H}})=1=[F]$, and that $C(T_{A/F})=[A]^{n/2}$
if $A$ has even degree (cf. Theorem 3). \hfill$\Box$
\vspace{0.1in}

\begin{rem}It follows from a comment given by Berhuy and Frings \cite[p.
4]{bf}, that if $A\in C.S(F)$ with odd degree $n$, then the second
trace form $T_{A/F}$ is Witt equivalent to $T_{M_n(F)}$ (the
second trace form for the matrix algebra of dimension $n$ over
$F$). As a consequence, the statements in Theorem 2 and Theorem 3
can be made more explicit in the case that the algebras $A_1,\
A_2\in C.S(F)$ both have odd degrees $n_1$ and $n_2$,
respectively. Indeed, we get upon invoking Proposition 1 that in
this case
\[
T_{A_1\otimes A_2/F}=
  \begin{cases}
    \frac{n_1^2n_2^2-1}{2}\mathbb{H} & \text{ if \ } n_1 n_2 \equiv 1,\ 7\mod 8, \\
   \left[1,1\right]+\frac{n_1^2n_2^2-3}{2}\mathbb{H} & \text{ if \ } n_1 n_2\equiv 3,\ 5 \mod
    8.
  \end{cases}\]
The corresponding invariants of $T_{A_1\otimes A_2/F}$ thus become
of the form
\[(Arf(T_{A_1\otimes A_2/F}), C(T_{A_1\otimes A_2/F}))=
\begin{cases}
(0, 1) & \text{ if \ } n_1 n_2 \equiv 1,\ 7 \mod 8, \\
(\left[1,1\right], ((1,1))_F) & \text{ if \  } n_1 n_2\equiv 3,\ 5
\mod 8.
  \end{cases}\]
\end{rem}

\end{document}